\documentclass[final, online]{aumbart}

\usepackage{amsmath}
\usepackage{amssymb}
\usepackage{amsthm}
\usepackage{enumerate}
\usepackage[mathscr]{eucal}
\usepackage[utf8]{inputenc}
\usepackage{graphicx}
\usepackage{url}
\usepackage{latexsym}

\newtheorem{thm}{Theorem}[section]
\newtheorem{lem}[thm]{Lemma}
\newtheorem{cor}[thm]{Corollary}
\newtheorem{prop}[thm]{Proposition}

\newtheorem{pr}[thm]{Problem}

\def\A{\bold A}

\def\D{\bold D}
\def\I{\bold I}
\def\M{\bold M}

\def\L{\bold L}
\def\R{\bold R}
\def\S{\bold S}
\def\X{\bold X}

\newcommand{\cat}[1]{\boldsymbol{\mathscr{#1}}}
\newcommand{\CA}{{\cat A}}
\newcommand{\CD}{{\cat D}}
\newcommand{\CCD}{{\cat D}}
\newcommand{\CF}{{\cat F}}

\newcommand{\CK}{{\cat K}}
\newcommand{\CL}{{\cat L}}
\newcommand{\CM}{{\cat M}}
\newcommand{\CR}{{\cat R}}
\newcommand{\CS}{{\cat S}}
\newcommand{\CX}{{\cat X}}

\def\twiddle#1{\smash{\underset{\smash{\raise.2ex\hbox{$\sim$}}}
              {#1}}\vphantom{#1}}
\newcommand{\thT}{\twiddle 3}
\newcommand{\twT}{{\twiddle 2}}
\newcommand{\MT}{{\twiddle M}}
\newcommand{\DT}{{\twiddle D}}
\newcommand{\RT}{{\twiddle R}}
\newcommand{\ZT}{{\twiddle Z}}
\newcommand{\Tp}{\mathscr{T}}
\newcommand{\T}{\mathscr{T}}
\newcommand{\x}{\text{\~x}}
\DeclareMathOperator{\fin}{fin} \DeclareMathOperator{\dom}{dom}
\DeclareMathOperator{\id}{id} \DeclareMathOperator{\graph}{graph}
\DeclareMathOperator{\Inv}{Inv} \DeclareMathOperator{\Pol}{Pol}
\DeclareMathOperator{\End}{End}
\DeclareMathOperator{\ISP}{\mathbb{ISP}}
\DeclareMathOperator{\IScP}{\mathbb{IS} _{\mathrm{c}}\mathbb{P}^+}

\begin{document}
\begin{frontmatter}
\title{On selected developments in the theory of natural dualities}
{}
\author{Miroslav Haviar}
{Faculty of Natural Sciences, M Bel University, \\ Tajovsk\'eho 40,
974~01 Bansk\'a Bystrica, Slovakia} {miroslav.haviar@umb.sk} {The
author gratefully acknowledges support from Slovak grant VEGA
1/0337/16.}

\keywords{Natural duality, entailment, endodualisability,
endoprimality, full and strong dualities}
\msc{Primary 08C20; 
Secondary 06D50}

\begin{abstract}
This is a survey on  selected developments in the theory of natural
dualities where the author had the opportunity to make with his
foreign colleagues several breakthroughs and move the theory
forward. It~is aimed as author's reflection on his works on the
natural dualities in Oxford and Melbourne over the period of twenty
years 1993-2012 (before his attention with the colleagues in
universal algebra and lattice theory has been fully focused on the
theory of canonical extensions and the theory of bilattices). It~is
also meant as a remainder that the main problems of the theory of
natural dualities, \emph{Dualisability Problem} and
\emph{Decidability Problem for Dualisability}, remain still open.

\emph{Theory of natural dualities} is a general theory for
quasi-varieties of algebras that generalizes `classical' dualities
such as \emph{Stone duality} for Boolean algebras, \emph{Pontryagin
duality} for abelian groups, \emph{Priestley duality} for
distributive lattices, and \emph{Hofmann-Mislove-Stralka duality}
for semilattices. We present a brief background of the theory and
then illustrate its applications on our study of Entailment Problem,
Problem of Endodualisability versus Endoprimality
and then a famous Full versus Strong Problem with related
developments.
\end{abstract}
\end{frontmatter}

\section{Introduction}
In 1936 M.H. Stone published a seminal work on duality theory,
exhibiting {a dual equivalence between the category of all Boolean
algebras and the category of all Boolean spaces}~\cite{St36}. Almost
at the same time L. Pontryagin showed that {the category of abelian
groups is dually equivalent to the category of compact topological
abelian groups} \cite{Po34a}, \cite{Po34b}. The most important step
toward the development of general duality theory was Priestley's
duality for distributive lattices: the category of all distributive
lattices was shown to be dually equivalent to the category of all
compact totally-order disconnected ordered topological spaces (since
then called {Priestley spaces})~\cite{Pr70}, \cite{Pr72}. Shortly
after that, K.H.~Hofmann, M.~Mislove and A.~Stralka developed a
duality for semilattices~\cite{HMS74}.
 The
general duality theory, called \emph{Natural duality theory}, grew
out from these four
 dualities,
in a monumental work by B.A. Davey and H. Werner~\cite{DW83}.  Its
rapid development over the next two decades is covered in the survey
papers by B.\,A.~Davey \cite{D93} and  by H.\,A.~Priestley
\cite{P93}, and in the monographs  by D.\,M.~Clark and B.\,A.~Davey
\cite{CD98} and by J.\,G.~Pitkethly and B.\,A.~Davey \cite{PiD05}.
The author's focus here is on selected developments in the theory
over the period of twenty years 1993-2012 where he had the
opportunity and privilege to make, mainly with H.\,A.~Priestley and
B.\,A.~Davey in Oxford and Melbourne, certain breakthroughs and move
the theory forward.

 The theory has proven to be a valuable tool in algebra, algebraic logic,
certain parts of computer science, and even in theoretical physics
as demonstrated by the author's survey in this journal on free
orthomodular lattices~\cite{H18}. This year's second (and expectedly
final) survey is also meant as a remainder that the main problems of
the theory, the \emph{Dualisability Problem} and the {Decidability
Problem for Dualisability}, remain still open.

Generally speaking, the theory of natural dualities concerns the
topological representation of algebras. The main idea of the theory
is that, given a quasi-variety $\CA = \ISP (\M)$ of algebras
generated by an algebra $\M$, one can often find a topological
relational structure $\MT$ on the underlying set $M$ of $\M$ such
that a dual equivalence exists between $\CA$ and a suitable category
$\CX$ of topological relational structures of the same type as
$\MT$. Requiring the relational structure of $\MT$ to be {\it
algebraic over $\M$}, all the requisite category theory ``runs
smoothly" (we refer to~\cite{CD98}). A uniform way of representing
each algebra $\bold A$ in the quasi-variety $\CA$ as an algebra of
continuous structure-preserving maps from a suitable structure
$\bold X\in \CX$ into $\MT$ can be obtained. In particular, the
representation is relatively simple and useful for free algebras in
$\CA$ as was demonstrated also in~\cite{H18}.

The motivation for the natural duality theory goes back to the
question ``Why in 1614 did the Scottish philosopher and
mathematician John Napier, Laird of Merchiston in Scotland, invent
the logarithm?" (\cite{Dcomm}). To quote from his 1619
book~\cite{N}:

{\it``Seeing there is nothing {\rm(}right well-beloved Students of
the Mathematics{\rm)} that is so troublesome to mathematical
practice, nor that doth more molest and hinder calculators, than the
multiplications, divisions, square and cubical extractions of great
numbers, which besides the tedious expense of time are for the most
part subject to many slippery errors, I began therefore to consider
in my mind by what certain and ready art I might remove those
hindrances. $\ldots$ I found at length some excellent brief rules
$\ldots$ which together with the hard and tedious multiplications,
divisions, and extractions of roots, doth also cast away from the
work itself even the very numbers themselves that are to be
multiplied, divided and resolved into roots, and putteth other
numbers in their place which perform as much as they can do, only by
addition and subtraction, division by two or division by three."\/}

A \emph{natural duality} is a form of logarithm which is applied to
algebraic structures rather than to numbers: it takes difficult
problems concerning algebras and converts them into simpler yet
equivalent problems concerning completely different mathematical
structures just as a logarithm converts a difficult multiplication
of positive real numbers into a simpler yet equivalent addition of
entirely different (and not necessarily positive) real numbers.
Given a finite algebra $\A$, a natural duality based on $\A$ is the
exact analogue of a logarithm, $\log_a$, to the base $a$ for
some positive real number $a \ne 1$ and $\A$ is said to 
\emph{admit a natural duality} if a natural duality based on $\A$
exists. Just as $\log_a$ does not exist if $a$ is not positive or $a
= 1$, a natural duality based on $\A$ need not exist. (\cite{Dcomm})

In Section 2 we present a brief background of the theory of natural
dualities with its main two open problems, the \emph{Dualisability
Problem} and the {Decidability Problem for Dualisability}. In
Sections 3 and 4 we illustrate the application of the theory on the
study of entailment and endodualisability developed by the author in
a close collaboration with H.A. Priestley and B.A. Davey. In Section
5 we give an overview of later developments of the theory in the
author's collaboration with B. Davey's research group, where our
focus is mainly on a famous Full versus Strong Problem.

\section{The basic scheme of the theory of natural
dualities and its main open problems}

We now  recall the basic scheme of the theory more precisely. Let
$\M = (M;F)$ be a finite algebra.
Let  $\MT = (M;G,H,R,\Tp)$ be a discrete topological structure,
i.e.~a non-empty set $M$ endowed with (finite) families $G$, $H$ and
$R$ of operations, partial operations and relations, respectively,
and with a discrete topology~$\Tp$.  We recall that the graph of an
$n$-ary (partial) operation $g\! : M^n\to M$ is the $(n+1)$-ary
relation
$$
\graph(g) = \{\,(x_1,\dots,x_n,g(x)) \mid (x_1,\dots,x_n)\in M^n\,\}
\subseteq M^{n+1}.
$$
We say that  the structure~$\MT$ is \emph{algebraic over\/}  $\M$ if
the relations in $R$ and the graphs of the operations and partial
operations in $G\cup H$ are subalgebras of appropriate powers
of~$\M$. Hence a unary (partial) operation is algebraic over $\bold
M$ if and only if it is a (partial) endomorphism of $\M$.

Let $\CA = \ISP (\M)$ be the quasi-variety generated by a finite
algebra $\M$ and assume that  $\MT = (M;G,H,R,\Tp)$ is algebraic
over $\M$. Let $\CX = \IScP(\MT)$ be the `topological quasi-variety'
generated by $\MT$, i.e.~the class of all structures which are
embeddable as closed substructures into powers of $\MT$. For any
algebra $\A \in \CA$, let $D(\A)$ denote the set of all
$\CA$-homomorphisms $\A\to \M$. Since $\MT$ is algebraic over $\M$,
$D(\A)$ can naturally be understood as a substructure of $\MT^A$,
and so as a member of $\CX$.

Let $X\subseteq M^I$ for some non-empty set $I$ and let $r\subseteq
M^n$ be an $n$-ary relation on $M$. We say that a map $\varphi\! :
X\to M$ preserves the relation $r$ if $ [\varphi
(\x_1),\dots,\varphi (\x_n)] \in r $ for all $\x_1=(x_{1i})_{i\in
I},\dots,\x_n=(x_{ni})_{i\in I}$ such that $[x_{1i},\dots,x_{ni}]\in
r$ for every $i\in I$. We say that $\varphi$ preserves an $n$-ary
(partial) operation if $\varphi$ preserves its graph as an
$(n+1)$-ary relation.

Let $\X$ be a structure in $\CX$. By an $\CX$-morphism $\varphi\! :
\X \to \MT$ we mean a continuous structure-preserving map, i.e.~a
continuous map preserving all (partial) operations in $G\cup H$ and
all relations in $R$. Let $E(\X)$ be the set of all $\CX$-morphisms
$\X \to \MT$. Again, since $\MT$ is algebraic over $\M$, $E(\X)$ can
be understood as a subalgebra of $\M^X$, i.e. a~member of $\CA$.

The (hom-)functors $D\! : \CA \to \CX$ and $E\! : \CX\to \CA$ are
contravariant and dually adjoint. Moreover, for any $\A\in \CA$ and
for any $\X \in \CX$, we have maps  $e_A\! : \A \to ED(\A)$ and
$\varepsilon_X\! : \X \to DE(\X)$ given by evaluation, {\it
viz.}
\begin{align} e_A(a)(h) &= h(a)\quad\text{ for every } a\in A \text{ and }
h\in D(\A)
 , \notag\\
\varepsilon_X(y)(\varphi) &= \varphi (y)\quad\text{ for every } y\in
X \text{ and }\varphi \in E(\X),\notag
\end{align}
which   are embeddings. We say 
that $\MT$ \emph{yields a pre-duality on} $\CA$.

Let $\MT = (M;G,H,R,\Tp)$ be an algebraic structure over $\M$, so
that $\MT$ yields a pre-duality on $\CA = \ISP (\M)$. We say that
\emph{$\MT$ yields a natural duality on $\CA$} if for every $\A\in
\CA$ the embedding $e_A$ is an isomorphism, i.e.~the evaluation maps
$e_A(a) \ (a\in A)$  are the only $\CX$-morphisms from $D(\A)$ to
$\MT$; we notice that they represent then the elements $a$ of $\A$.
Sometimes we say that \emph{$G\cup H\cup R$ yields a
{\rm(}natural\/{\rm)} duality on\/} $\CA$ or that \emph{$\MT$ is
dualisable}. We further say that $\MT$ (or $G\cup H\cup R$) yields a
\emph{full duality\/} on $\CA$ if $\MT$ yields a duality on $\CA$
and for every $\X\in \CX$ the embedding $\varepsilon_X$ is also an
isomorphism. In such a case the categories $\CA$ and $\CX$ are
dually equivalent via categorical anti-isomorphisms $D$ and $E$
which are inverse to each other. Finally, we say that $\MT$ (or
$G\cup H\cup R$) yields a \emph{strong duality\/} on $\CA$ if $\MT$
is injective in the category~$\CX$ (with respect to embeddings). A
famous \emph{Full versus Strong Problem}, which dated back to the
beginnings of the theory of natural dualities and was open for about
twenty-five years asked:

\begin{pr}  \rm{(}{Full versus Strong Problem}\rm{)}
{\it Is every full duality strong?}
\end{pr}

We have not claimed above that it is always possible, for a given
algebra $\M$, to choose a  structure $\MT$ on $M$ yielding a
 duality on $\ISP (\M)$. In fact, the main
problem of the theory of natural dualities, the \emph{Dualisability
Problem}, remains still open:

\begin{pr}  \rm{(}{Dualisability Problem}\rm{)}
{\it Which finite algebras are dualisable?}
\end{pr}

At present, the Dualisability Problem seems to be unsolvable
(cf.~\cite[page viii]{PiD05}). There are algebras $\M$ which fail to
be dualizable (we refer to~\cite{DW83} or \cite{D93}). However, for
a very wide range of algebras dualities do exist.  For example,
 the NU-Duality  Theorem (\cite{DW83}, Theorem 1.18 or \cite{D93}, Theorem 2.8) guarantees that a duality
on $\ISP (\M)$ is available whenever $\M$ has a lattice reduct. Many
further  theorems which say how to choose an appropriate structure
$\MT$ on $M$ to obtain a  duality, or a strong (thus full) duality,
on $\ISP (\M)$ can be found in~\cite{CD98} and in~\cite{PiD05}. The
Dualisability Problem might be formally undecidable, and in fact,
the ``holy grail" (cf.~\cite[page viii]{PiD05}) of some
natural-duality theoreticians is the {Decidability Problem for
Dualisability}:

\begin{pr}  \rm{(}{Decidability Problem for Dualisability}\rm{)}
{\it Is there an algorithm for deciding whether or not any given
finite algebra is dualisable?}
\end{pr}

\section{Entailment in natural dualities and our solution of the Entailment problem}

Again assume a structure $\MT = (M;G,H,R,\Tp )$ is algebraic over a
finite algebra $\M$ and let $r$ be an $n$-ary algebraic relation on
$M$ (i.e.~a subalgebra of $\M^n$). We say that the structure $\MT$,
or more often just $G\cup H\cup R$, \emph{entails\/} $r$ if for
every $\X\in \CX$, each $\CX$-morphism $\varphi\! : \X\to \MT$
preserves~$r$; we write $G\cup H \cup R \vdash r$. For relations $r$
and $s$ we write $r \vdash s$ in place of $\{ r\} \vdash s$.  We say
that $G\cup H\cup R$ entails an $n$-ary (partial) operation $h$ if
it entails its graph  as an $(n+1)$-ary relation, and that it
entails a set $R'$ of relations and (partial) operations if it
entails each $r\in R'$.

\subsection{Test Algebra Lemma and the Entailment problem}

Central to the identification of the relations entailed from certain
set $G \cup H \cup R$ is so-called Test Algebra Lemma. (It is
formulated in entailment terms in \cite{DP96}, Lemma 2.3 and in
\cite{CD98}, Lemma 8.1.3.) We present this statement and we notice
that $\bold s$ always denotes the algebraic relation $s$ considered
as an algebra in $\CA$.

\begin{thm} \rm{(}{{Test Algebra Lemma}}\rm{)} Let $\M$ be a finite algebra,
let $G$, $H$, $R$ be, respectively, sets (possibly empty) of
operations, partial operations and relations which are algebraic
over~$\M$, and let $s$ be an algebraic relation. Then the following
are equivalent:
\begin{itemize}
\item[(1)]  $G \cup H \cup R$ entails $s$;
\item[(2)] $G \cup H \cup R$ entails $s$ on $D(\bold s)$.
\end{itemize}
Moreover, $G\cup H\cup R$ entails $s$ whenever $G\cup H\cup R$
yields a duality on $\bold s$.
\end{thm}

We often use the term \emph{test algebra} for an algebra $\bold A\in
\ISP(\M)$ witnessing the  failure of the  structure $\MT$ to yield a
duality on $\ISP(\M)$.

It is important that provided a set $G\cup H\cup R$ yields a duality
on $\CA$ then the duality is not destroyed by deleting from $G \cup
H \cup R$ any element which is entailed by the remaining members.
This is the key to obtaining so-called  \emph{economical dualities}
which are easy to work with. A full discussion of the central role
played by entailment in duality theory is presented in the paper
\cite{DHP:syntax}. In this paper we solved the Entailment Problem of
duality theory that was formulated as follows:

\begin{pr} \rm{(}{Entailment Problem}\rm{)}
{\it Find an intrinsic description of the relations entailed by $G \cup H \cup R$.}
\end{pr}

This problem was formulated as the first open problem of the natural
dualities in the famous survey paper~\cite{D93}. When this problem
was firstly introduced, it was expected that the solution would be a
semantic one in terms of a preservation theorem providing a list of
finitary constructs which preserve entailment. By this is meant that
if $(G \cup H \cup R) \vdash s$ then $s$ would be obtainable from
the set $G \cup H \cup R$ via a finite sequence of finitary
constructs. In our solution to the problem in \cite{DHP:syntax} we
indeed firstly used a semantic approach, which
was similar to the characterisation of the well-known \emph{clone
closure} $\Inv(\Pol($R$))$ of a set of relations $R$ (all
`invariants' of `polymorphisms' preserving~$R$) originally obtained
in the famous pair of papers~\cite{BKKR69} by V.~Bodnar{\v c}uk,
L.A.~Kalu{\v z}nin, V.N.~Kotov and B.A.~Romov. Later on, we noticed
that our semantic solution also arises as a direct application of a
syntactic solution: a description  of relations entailed by $G \cup
H \cup R$ in terms of the first-order formul{\ae} of the language
with equality, $\CL_{\MT}$, associated with $\MT$. An important step
towards  the solution
 was  the recognition
that on a given set $\Omega $ of finitary algebraic relations on
$\M$ the map $ R \longmapsto \overline{R} := \{ \, s \in \Omega \mid
R \vdash s \,\} $ is a closure operator (\emph{entailment closure}).
And also the recognition that this closure operator is algebraic, in
the sense that the closure of any set $R$ is the union of the
closures of its finite subsets (so that the lattice  of closed sets
is algebraic).
 This provided indirect evidence for
a positive solution to the Entailment Problem.

\subsection{Our syntactic solution of the Entailment problem}

In \cite{DP96} the important fact that entailment closure is
algebraic was deduced as a corollary of the Test Algebra Lemma. In
the paper \cite{DHP:syntax} we extended the Test Algebra Lemma,
upgrading it to the Test Algebra Theorem. This theorem provides our
syntactic solution to the Entailment Problem:

\begin{thm} {\rm(}The Test Algebra Theorem or Entailment in the duality sense{\rm)}
Let $\M$ be a finite algebra and let a structure $\MT = (M;
G,H,R,\Tp )$ be algebraic over $\M$.  Then the following are
equivalent:
\begin{itemize}
\item[(1)]  $G \cup H \cup R$ entails $s$;
\item[(2)] $G \cup H \cup R$ entails $s$ on $D(\bold s)$;
\item[(3)] some finite subset of $G \cup H \cup R$ entails $s$ on $D(\bold s)$;
\item[(4)] $s = \{\, (u(\rho _1), \dots, u(\rho _n)) \mid
u : D(\bold s) \to M \text{ preserves } G \cup H \cup R \,\}$;
\item[(5)] there exists a primitive positive formula
$\Phi(x_1,\dots,x_n)$ in the language $\cat L_{\MT}$ such
that
{\begin{itemize}
\item[(i)] $D(\bold s) \vdash \Phi (\rho_1,\dots,\rho_n)$ and

\item[(ii)] $s = \{\,(c_1,\dots,c_n) \in M^n \mid M \vdash \Phi(c_1,\dots,c_n)
\,\}$.
\end{itemize}}
\end{itemize}
\end{thm}

The most important part of our  syntactic solution is that $(G \cup
H \cup R ) \vdash s$ if and only if there is a primitive positive
formula $\Phi $ in the language $\CL_{\MT}$ such that $s$ may be
obtained from $G\cup H\cup R$ via a \emph{primitive positive
construct}. We may take $\Phi$ to be the primitive positive type of
$\rho_1, \dots, \rho_n$ in $D(\bold s)$.

In duality theory, a set $R$ of finitary algebraic relations on a
finite algebra $\M$ entails a finitary algebraic relation $s$ on the
powers of $\MT$ (which are the duals of free algebras in the
associated quasivariety $\CA$; see, for example, \cite{DW83}) if and
only if $s$ can be obtained from $R$ in the clone-theoretic case.

Therefore applying our results in the clone setting we derive a
famous consequence due to V.~Bodnar{\v c}uk, L.A.~Kalu{\v z}nin,
V.N.~Kotov and B.A.~Romov~\cite{BKKR69}:

\begin{thm} {\rm(}Entailment in the clone sense{\rm)}
Let $R$ be a family of finitary relations on a finite set $M$
and let $s \subseteq M^n$.
Then the following are equivalent:
\begin{itemize}
\item[(1)] $s \in  \text{\rm Inv(Pol($R$))} $;
\item[(2)] $R$ entails $s$ on $M^s$;
\item[(3)] $s = \{( u(\rho _1), \dots, u(\rho _n)) \mid u: M^s \to M \text{ preserves } R\}$;
\item[(4)] there is some finite structure $Z$ of type $(M;R)$ and elements $z_1, \dots, z_n \in Z$ such that
$s = \{\, (u(z _1), \dots, u(z _n)) \mid u: Z \to M  \text{
preserves } R\,\}$;
\item[(5)]  $s = \{\,(c_1,\dots,c n) \in M^n \mid M \vdash \Phi(c_1,\dots,c_n)\,\}$ for some primitive positive formula
$\Phi(x_1,\dots,x_n)$ {\rm(}in the language of the relational structure $(M;R)${\rm)}.
\end{itemize}
\end{thm}

\subsection{Our semantic solution of the Entailment problem}

Through the Test Algebra Theorem we are able to  convert our
syntactic solution to the Entailment Problem to a semantic solution,
so obtaining a set of constructs sufficient to describe entailment.
We only summarise the results below and sketch the main steps of our
semantic solution while for all details of it and definitions of the
constructs we refer to our paper~\cite{DHP:syntax} or to~\cite[2.4.5
and 9.2.1]{CD98}.

In case $G \cup H = \emptyset$, the list of entailment constructs
may be taken to be: \emph{trivial relations}, \emph{repetition
removal}, \emph{intersection}, \emph{product}, and \emph{retractive
projection} (in which the natural projection map is required to be a
retraction). As a consequence in the clone setting we have the
result of~\cite{BKKR69} that $\text{Inv(Pol($R$))}$ can be obtained
from $R$ by a finite number of applications  of trivial relations,
intersection, repetition removal, product and projection.

As is well known, arbitrary projection is not necessarily an
allowable construct on structures of the form $D(\bold A) ={\CA}
(\bold A, \M)$. If it were,  we could form  the relational product
of two relations, which is not guaranteed to lift to  structures
$D(\bold A)$ which are not full powers.  This explains why a set $R$
of algebraic relations on $\M$ which determines the clone of term
functions on $\M$ will not necessarily yield a duality on~$\CA$.
This is illustrated in \cite[p.102]{D93} in case $\CA$ is the
variety $\CK$ of Kleene algebras; for a more extended discussion we
refer to~\cite[Section 5]{DP96} or \cite{DHP:kleene}.

Our semantic solution to the Entailment Problem in~\cite{DHP:syntax}
was carried out in two stages. Firstly, we showed that the second
dual $ED(\bold s)$ of an algebraic relation $s$ can be
\emph{concretely} constructed from $G \cup H \cup R$, whether or not
$G \cup H \cup R$ entails $s$ (for details again
see~\cite{DHP:syntax} or~\cite[2.4.5 and 9.2.1]{CD98}). Secondly, we
showed that if $G \cup H \cup R$ entails $s$ then $s$ can be
obtained from this second dual $ED(\bold s)$ by a {retractive
projection}, which is a {bijective projection} in case $G \cup H
\cup R$ yields a duality on $\bold s$.

To explain the latter concepts, given an $m$-ary algebraic relation
$r$ on $M$ and an injective mapping
 $\eta: \{1,\dots,n\} \to \{1,\dots,m\}$ ($n \leq m$) we define the relation
$$
r_\eta = \{\,(c_1,\dots,c_n) \in M^n \mid (\exists d_1 \dots d_m\in
M)\ (d_1,\dots,d_m) \in r \text{ and } c_i = d_{\eta (i)} \ (1 \leq
i \leq n) \, \}
$$
(it can be alternatively denoted as the projection $P_{\eta (1),
\dots, \eta (n)}(r)$ of $r$ into its coordinates $\eta (1), \dots,
\eta (n)$). Then we say that the relation $s:=r_\eta$ is a
\emph{retractive projection of $r$} if the natural projection map
$p:  \bold r \to \bold s $ is a retraction, that is, there is a
homomorphism $q: \bold s \to \bold r$ such that $p\circ q = \id_s$.
It is called a \emph{bijective projection} (as introduced by
L.~Zadori \cite{Z}) if moreover $q \circ p = \id_r$. A retractive
projection derived from an injection of $\{1, \dots, m-1\}$ into
$\{1, \dots, m\}$ is called a \emph{$1$-step retractive projection
of $r$}.

Consider $G$, $H$ and $R$ as before and let now
$Z = \{z_1,\dots,z_k\}$ be a finite substructure of $M^T$, for
some non-empty set $T$.
By the \emph{graph of $E(Z)$} (with respect to $G \cup H \cup R$)
we mean the relation
$$
G[E(Z)] := \{\,(u(z_1),\dots,u(z_k)) \in M^k \mid u: Z \to
M \text{ preserves } G \cup H \cup  R \,\}.
$$
Thus the graph of $E(Z)$ is simply  $E(Z)$, given a fixed labelling
of $Z$. We showed that if $Z$ is a finite  subset of $M^T$ for some
non-empty set $T$ which is \emph{hom-closed} (for details
see~\cite[p.~66]{CD98}), then the relation $G[E(Z)]$ can be
concretely constructed from $G \cup H \cup R$.

For an $n$-ary algebraic relation $s$ we take $Z := D(\bold s)$ to
be the dual of the algebra $\bold s$ and enumerate  its elements as
$\{\rho _1,\dots ,\rho _n,\Tp _1,\dots,\Tp _m\}$.  We then assume
that
$$
G[\bold s] := \{\,(\rho_1(a),\dots,\rho_n(a),\Tp _1(a),\dots,\Tp
_m(a)) \in M^{n+m} \mid a\in s\,\}
$$
 encode the evaluation maps from $D(\bold s)$ to $M$.  It is evident that $G[\bold s]$ is in bijective
correspondence with  $s$ itself. Now we have that if $G \cup H \cup
R$ yields a duality on $\bold s$ then $G[ED(\bold s)]$ necessarily
coincides with $G[\bold s]$.  It is helpful to employ the intuition
that the relation $G[ED(\bold s)] \setminus G[\bold s]$ can be
thought of as a measure of how far $G \cup H \cup R$ is from
yielding a duality on $D(\bold s)$.

Since by the Test Algebra Theorem we have that an algebraic relation
$s$ is the  retractive projection of $G[ED(\bold s)]$ onto its first
$n$ coordinates, where the dual $D(\bold s)$ of $s$ is labelled as
above, we immediately have:

\begin{lem} Let $\bold s \leq \bold M^n$
and $G \cup H \cup R$ entail $s$.
 Then $s$ is a retractive
projection of the graph   $G[ED(\bold s)]$ of $ED(\bold s)$.
 \end{lem}

A number of consequences can be deduced.  The first is the desired
Semantic Entailment Theorem of~\cite{DHP:syntax}:

\begin{thm} {\rm(}Semantic Entailment
Theorem{\rm)}
Let $R$ be a set of algebraic relations on a finite
set $M$, let $s$ be an algebraic relation on $M$  and let  $R \vdash
s$.  Then $s$ can be obtained from $R$ by a finite number of
applications of product, intersection, trivial relations and
repetition removal, followed by one application of retractive
projection.
\end{thm}

If a set $R$ of algebraic relations on a finite set $M$ is such that
 $R \vdash s$ for every  algebraic relation $s$ on $M$, then we say that
$R$ is \emph{entailment-dense}. The following result, that can be
derived from our semantic solution, was (independently to our
investigations) discovered by L.~Z\'adori~\cite{Z}:

\begin{thm} {\rm(}Special Semantic Entailment
Theorem{\rm)}
 Let $R$ be a set of algebraic relations on a finite
set $M$ and let $s$ be an algebraic relation on $M$.
\begin{itemize}
\item[(a)]  If $R$ yields a duality on $\bold s$, then $s$ can be
 constructed from $R$ by a finite number of applications of product,
intersection, trivial relations, repetition removal and bijective
projection.
\item[(b)]  The following are equivalent:
{\begin{itemize}
\item[(i)] $R$ yields a duality on every finite algebra in $\CA$;
\item[(ii)] $R$ is entailment-dense;
\item[(iii)] every algebraic relation $s$ on $M$ can be constructed from
$R$ by a finite number of applications of product, intersection,
trivial relations, repetition removal and bijective projection.
\end{itemize}}
\end{itemize}
\end{thm}

\section{Endoprimality and endodualisability in theory and practice}

The relationship between duality entailment and clone-entailment is
rather complex. It is known  that it is possible for $G\cup H\cup R$
to clone-entail every finite algebraic relation on $\M$ but to fail
to dualise $\M$, but the circumstances under which this phenomenon
occurs, and what it signifies, are still obscure. In particular, we
may ask what it means for $\M$ to be \emph{endoprimal} but not
\emph{endodualisable} (we refer to definitions of these concepts
below). More explicitly, we may ask what it means for some finitary
algebraic relation $r$ on $\M$ to be clone-entailed but not entailed
by (the graphs of) the endomorphisms of $\M$. From a semantic
viewpoint, a clear difference can be seen: clone-entailment allows
all relational products, whereas duality entailment allows only
\emph{homomorphic relational products} (for details
see~\cite{DHP:syntax} or~\cite[9.2.1]{CD98}). Thus one may expect
relational products appearing in the construction of $r$ from the
endomorphisms of $\M$ to be non-homomorphic relational products.
Exactly how this behaviour happens in general is not clear.

\subsection{Endoprimality versus endodualisability}

In \cite{DHP:kleene} we showed  that the relationship between the
two entailment concepts also lies at the heart of the relationship
between endoprimality and endodualisability. This was nicely
demonstrated by the Kleene algebra examples. We note that Kleene
algebras were already known to illustrate the distinction between
entailment in the clone sense and in the duality sense - we refer to
\cite[p.~87]{D93}, \cite[Section 5]{DP96} and~\cite[pp.
272--273]{CD98}. In \cite{DHP:kleene} we give a complete description
of endodualisable and endoprimal finite Kleene algebras from the
quasi-variety $\ISP(\bold 4)$ and show that there is a plentiful
supply of finite Kleene algebras which are endoprimal but not
endodualisable.

Let $\M = (M;F)$ be any algebra. The  algebra $\M$ is called
$k$-\emph{ endoprimal} ($k\geq 1$) if every $k$-ary $\End(\bold
M$)-preserving function on $\M$ is a term function of $\M$. Algebras
which are $k$-endoprimal for every $k\geq 1$ are called
\emph{endoprimal}. A finite algebra $\M$ is \emph{endodualisable} if
$\End(\M)$ yields a duality on the quasivariety $\ISP(\M)$.

The relationship between endodualisability on one hand, and
endoprimality and $k$-endoprimality on the other hand, has been
explored, successively, in \cite{DHP:endo}, \cite{D96},
\cite{DPi97}, \cite{HP99a} and~\cite{DHP:kleene}. It has been shown
that in many quasivarieties a finite algebra is endoprimal if and
only if it is endodualisable (we refer to \cite{DPi97}, \cite{HP99b}
and the papers cited therein).

In \cite{DHP:endo} we started an intensive study of a general
relationship between endodualisability and endoprimality by the
following result:

\begin{thm} {\rm(}Endoprimality versus endodualisability for distributive lattices{\rm)
}
 Let $\L = (L;\lor,\land)$ be a finite non-trivial  distributive lattice. The following are
equivalent:
\begin{itemize}
\item[(1)] $\L$ is {\rm 3}-endoprimal;
\item[(2)] $\L$ is endoprimal;
\item[(3)] $\L$ is endodualisable;
\item[(4)] the retractions of $\L$ onto $\{0,1\}$ together with the constants $0, 1$ yield a duality on $\ISP(\L)$;
\item[(5)]  $\L$ is not a Boolean lattice.
\end{itemize}
\end{thm}

In case of bounded distributive lattices we obtained a similar
result, the only difference is in Condition (1):

\begin{thm} {\rm(}Endoprimality vs endodualisability for bounded distributive lattices{\rm)}
 Let $\L = (L;\lor,\land,0,1)$ be a finite non-trivial bounded distributive lattice. The following are
equivalent:
\begin{itemize}
\item[(1)] $\L$ is {\rm 1}-endoprimal;
\item[(2)] $\L$ is endoprimal;
\item[(3)] $\L$ is endodualisable;
\item[(4)] the retractions of $\L$ onto $\{0,1\}$ together with the constants $0, 1$ yield a duality on $\ISP(\L)$;
\item[(5)]  $\L$ is not a Boolean lattice.
\end{itemize}
\end{thm}

The first examples of finite algebras which are endoprimal but not
endodualisable were found by B.A.~Davey and J.G.~Pitkethly in their
paper~\cite{DPi97}, among algebras with a semilattice reduct. Many
other such examples have been found  among Kleene algebras in our
paper~\cite{DHP:kleene}.

\subsection{A criterion for a finite endoprimal algebra to be endodualisable}

In the  paper \cite{HP99a} the  strategy for finding   endoprimal
algebras due to  B.A.~Davey and J.G.~Pitkethly~\cite{DPi97} is
further explored  in the finite case. A new theoretical tool, called
the \emph{Retraction Test Algebra Lemma}, is used to show that, in
many quasivarieties, endoprimality is equivalent to
endodualisability for finite algebras which are suitably related to
finitely generated free algebras. The main result of \cite{HP99a} is
the following theorem.

\begin{thm}{\rm(}Retraction Test Algebra Lemma{\rm)}
Let a finite algebra $\D$ be dualisable via the structure
$$
\DT =(D;\End(\D),s_1,\dots,s_m, \Tp)
$$
where $m\ge 1$ and $s_1,\dots,s_m$ are finitary algebraic relations
on $\D$. Let the algebras $\bold s_1,\dots,\bold s_m$ be retracts of
the $k$-generated free algebra $\bold F_{\CD}(k) \in \CD$ where
$\CD=\ISP (\D)$.

Then for any finite algebra $\M \in \CD$ which has $\D$ as a
retract the following are equivalent:
\begin{itemize}
\item[(1)] $\M$ is endoprimal;
\item[(2)] $\M$ is $k$-endoprimal;
\item[(3)] $\M$ is endodualisable.
\end{itemize}
\end{thm}

The result can  be applied to the (quasi-)varieties of distributive
lattices (with $k~=~3$), bounded distributive lattices ($k=1$),
finite vector spaces of dimension greater than one ($k=2$), Stone
algebras ($k=2$), abelian groups ($k=2$), sets ($k=3$),
 semilattices ($k=3$), lower-bounded semilattices ($k=2$) and
 median algebras ($k=3$), which have not been considered before
as regards endoprimality.

We explain the applications of our theorem above in several selected cases:

\smallskip

\noindent
{\bf Distributive lattices}

The class  $\CD$ of distributive lattices is the quasi-variety $\ISP
(\bold 2)$ generated by the $2$-element lattice $\bold 2 =
(\{0,1\};\lor,\land)$. It is well-known (by \emph{Priestley duality}
presented in~\cite{Pr70}, \cite{Pr72}) that  $\bold 2$ is dualisable
via the structure $\twT=(\{0,1\},0,1,\le ,\Tp)$ where  $\le$ is the
usual order on $\{0,1\}$. It is said that $\bold 2$ is \emph{almost
endodualisable} with $\le $ as the extra relation to the
endomorphisms in the dualising structure. We notice that $\le$ is,
as a distributive lattice, isomorphic to the $3$-element chain
$\bold 3$.

It is easy to check that the free algebras $\bold F_{\CD}(1) \cong \bold 1$ and $\bold F_{\CD}(2) \cong \bold
2^2$  do not have $\bold 3$ as a retract while the free algebra
$\bold F_{\CD}(3)$ does have $\bold 3$ as a retract.
All non-trivial distributive lattices $\bold L\in \CD$ have evidently $\bold 2$ as their retracts.
From our theorem above it therefore follows that a finite non-trivial distributive lattice $\L$
is endoprimal iff $\L$ is $3$-endoprimal iff $L$ is endodualisable.

\smallskip

\noindent
{\bf Stone algebras}

The class of Stone algebras is the quasi-variety $\ISP (\bold 3)$
generated by the $3$-element Stone algebra $\bold 3 =
(\{0,a,1\};\lor,\land,^\star,0,1)$ where $\{0,a,1\}$ is the
$3$-element chain and $0^\star =1$ and $a^\star =1^\star = 0$. It is
well known  that the structure $\thT = (\{0,a,1\},d,\preccurlyeq,\Tp
)$ yields a duality on the variety of Stone algebras (cf. e.g.
\cite[p.~105]{CD98}) where  $\preccurlyeq$ is the order
$\{(0,0),(a,a),(1,1),(1,a)\}$ and $\graph(d)=\{(0,0),(1,1),(a,1)\}$.
It means that $\bold 3$ is almost endodualisable with the extra
relation $\preccurlyeq$ which is isomorphic to the $4$-element chain
algebra $\bold 4$ in $\CS$. Now the smallest $k$-generated free
algebra in $\CS$
 having  $\bold 4$ as a retract is known to be $\bold F_{\CS}(2)$.
 Our theorem  can now be applied to Stone algebras having $\bold 3$
as a retract. The only Stone algebras which do not have $\bold 3$ as
a retract
 are the Boolean algebras (and these  are
endodualisable). It follows that a finite non-Boolean
Stone algebra $\L$ is endoprimal iff $\L$ is
$2$-endoprimal iff $\L$ is endodualisable.

\smallskip

\noindent
{\bf Median algebras}

The class of median algebras is  the quasi-variety $\CM=\ISP (\bold
M)$ generated by the $2$-element median algebra $\M= (\{0,1\};m)$ in
which  the ternary (median) operation $m$ satisfies the equations
$$
m(x,y,z)=m(y,x,z)=m(y,z,x), \ m(x,x,y)=x
$$
and
$$
m(m(x,y,z),u,v)=m(x,m(y,u,v),m(z,u,v)).
$$

The duality for  $\CM$ is given by the structure
$\MT=(\{0,1\};^*,0,1,\le,\Tp)$,
 where $^*$ is the automorphism reversing $0$ and $1$ and $\le$ is the usual order on $\{0,1\}$ (we refer, for example, to~\cite[p.~103]{CD98}).
 It follows that
 $\M$ is almost endodualisable with the extra relation $\le$ which can be considered as a
median algebra, say $\bold s$. In our paper~\cite{HP99a} we present a verification in terms of natural duals of the fact
that the smallest $k$-generated free algebra in $\CM$ which
has the algebra $\bold s$ as a retract is $\bold F_{\CM}(3)$.
Because any non-trivial median algebra $\L\in \CM$ has
$\bold M$ as a retract it immediately follows from
our theorem that a finite non-trivial median algebra $\L\in \CM$ is
endoprimal iff $\L$ is $3$-endoprimal iff $\L$ is
endodualisable.

\smallskip

\noindent
{\bf Abelian groups}

Our method allows us to identify also the finite endoprimal abelian
groups.  Starting from a finite abelian group $\A$, one can choose
$\CD$ and the generator $\bold D$ of $\CD$ in such a way that $\A
\in \CD$ and $\bold D$ is a retract of $\A$.  This enables us to
apply our theorem.

It is well-known that for any finite abelian group $\A$ there is a
cyclic group $\bold Z_m$ such that $\bold A\in \CA_m$ where
$\CA_m=\ISP(\bold Z_m)$ and $\bold Z_m$ is a direct factor, and
hence a retract, of $\A$. It was shown in~\cite{DW83} (we also refer
to \cite[p.~114]{CD98}) that the structure $\ZT_m = (Z_m;+,^-,0,\Tp
)$ yields a duality on the quasi-variety $\CA_m$. This means that
$\bold Z_m$ is almost endodualisable with $\graph(+)$ as the extra
relation, which is, as an algebra, isomorphic to $\bold Z_m^2$. We
have $\bold F_{\CA_m}(2) \cong \bold Z_m^2$. Hence for the finite
abelian group $\A$ and the associated quasivariety $\CA_m
=\ISP(\bold Z_m)$ we could apply  our theorem with $k=2$. It follows
that a finite abelian group $\bold A$ is endoprimal iff it is
$2$-endoprimal iff it is endodualisable.

\subsection{Endodualisable and endoprimal finite double Stone
algebras}

In the  paper \cite{HP99b}
we give a complete characterisation of the endoprimal  finite double
Stone algebras. In particular, we have
 shown that all of these algebras are endodualisable,
and found in every case the minimum value of $k$ for which
$k$-endoprimality forces endoprimality. Much more work was involved
in completing this  analysis than that for the other examples
considered in the paper \cite{HP99a}, and
further duality techniques were required.

Let us present a brief outline of the results. An algebra $\L =
( L;\lor,\land,^\star ,^+ ,0,1)$ is called a \emph{double Stone
algebra} if $(L;\lor,\land,^\star ,0,1)$ and $(L;\land,\lor,^+
,1,0)$ are Stone algebras. The  double Stone algebras form a variety
$\CD\CS = \ISP (\bold 4)$  which is generated by the $4$-element
chain algebra $ \bold 4 =  ( \{0,a,b,1\};\lor,\land,^\star ,^+ ,0,1)
$ where $0<a<b<1$ and
$$
1^\star=b^\star=a^\star=0,\ 0^\star = 1,\ 0^+=a^+=b^+=1,\ 1^+ = 0.
$$
The proper non-trivial subvarieties of $\CD\CS $ are generated
by the subdirectly irreducible subalgebras $\bold 2 = \{ 0,1\}$  and
 $\bold 3 = \{ 0,a,1\}$.
The variety $\ISP(\bold 2)$  is just the class of Boolean algebras,
while $\ISP(\bold 3)$ is the variety of regular double Stone
algebras, {\it alias\/} three-valued
Lukasiewicz algebras. 
An algebra is \emph{proper} precisely when it  has $\bold 4$ as a
retract. We have to consider  separately the algebras in $\ISP(\bold
4) \setminus \ISP(\bold 3)$, which we call proper double Stone
algebras, and algebras in $\ISP(\bold 3)$. Also, a further splitting
into cases is necessary, into algebras with non-empty core and
algebras with empty core. The \emph{core} of an algebra $\bold L$ in
$\CD\CS$ is  defined to  be
 $K(\bold L) =  \{\,x\in L\mid x^{\star} =0, \, x^+ =
1\,\}$. A finite algebra $\bold L$ has empty core if and only if
$\bold L$  has $\bold 2$ as a direct factor.  It is easily shown
that this occurs if and only if $\bold L \in \ISP(\bold 4 \times
\bold 2)$. Every $k$-generated free algebra $\bold F_{\CD\CS}(k)$ lies in the  subquasivariety  $\ISP(\bold 4 \times \bold 2)$.

The finite non-Boolean algebras in the variety $\ISP(\bold 3)$ are
exactly those of the form $\bold 3^m \times \bold 2^\ell$ ($m\geq
1$, $\ell \geq 0$). We can set up  a duality for $\ISP(\bold 3
\times \bold 2)$ in which the only
non-endomorphism is isomorphic to $\bold 3\times \bold 2^2$.

We proved the following result:

\begin{thm} {\rm(}Endodualisable finite double Stone
algebras{\rm)}
 Let $\bold L$ be a finite non-trivial double Stone algebra and express  $\bold L$ as $\bold J
\times \bold 2^\ell$  where $\bold J$ does not have $\bold 2$ as a
factor and $\ell \geq 0$.

Then $\bold L$ is endodualisable when $\bold L$ takes one of the
forms described below.
\begin{itemize}
\item[(1)] $\L$ has  non-empty core and $\L$ satisfies  the following equivalent conditions:
{\begin{itemize}
\item[(i)] $\L$ has $\bold 5$ as a retract;
\item[(ii)] $K(\L)$ is a non-Boolean lattice.
\end{itemize}}
\item[(2)] $\L$ is proper, $\bold J$ has $\bold 5$ as a retract and $\ell  \geq 2$.
\item[(3)] $\L$ is not proper and takes the form $\bold 3^m \times \bold 2^\ell$ where $m \geq 1$ and $\ell \geq 2$.
\item[(4)]  $\L$ is Boolean.
\end{itemize}
\end{thm}

Let $\bold L$ be a finite non-trivial and non-Boolean double Stone
algebra which is not shown by above theorem to be endodualisable and
assume that $\bold L$ is expressed as $\bold J \times \bold 2^\ell$
where $\bold J$ does not have $\bold 2$ as a factor. The following
cases arise:

\begin{itemize}
\item[(A)] $\bold L$ is a Post algebra of order $3$ (that is, $\bold L$ is not proper  and $\ell= 0$);
\item[(B)] $\bold L$ has a single factor $\bold 2$ (that is,
$\ell =1)$;
\item[(C)] $\bold L$ is proper, $K(\bold L) \ne \emptyset$ (that is, $\ell =0)$, and $\bold J$ does not have
$\bold 5$  as a retract;
\item[(D)]  $\bold L$ is proper,  $\bold J$ does not have $\bold 5$ as  a retract and $\ell\geq 2$.
\end{itemize}

We showed that $\bold L$ is not endodualisable in each of cases
(A)--(D), treating these in turn.

\begin{prop} {\rm(}Non-endodualisable finite double Stone
algebras, Case A{\rm)} Let $\bold L$ be a  finite Post algebra of
order $3$. Then
\begin{itemize}
\item[(1)] $\bold L$ is not endodualisable, with $\bold 2$ serving as a test
algebra;
\item[(2)] $\bold L$ is not $1$-endoprimal.
\end{itemize}
\end{prop}

\begin{prop} {\rm(}Non-endodualisable finite double Stone
algebras, Case B{\rm)} Let $\bold L=\bold J \times \bold 2 $ be a
finite non-Boolean double Stone algebra with exactly one factor
$\bold 2$. Then \begin{itemize}
\item[(1)] $\bold L$ is not endodualisable, with  $\bold 2^2$
serving as a test algebra;
\item[(2)] $\bold L$ is not $1$-endoprimal.
\end{itemize}
\end{prop}

 For case~(C) we showed that
 the algebra $\bold L$
is the retract of a power of a finite indecomposable algebra which
is not $3$-endoprimal.

\begin{prop} {\rm(}Non-endodualisable finite double Stone
algebras, Case C{\rm)}
 Let $\bold L$ be a finite proper double Stone
algebra with a non-empty core $K(\bold L)=[a,b]$ {\rm(}$a<b${\rm)}
which is a Boolean lattice. Then $\bold L$ is not $3$-endoprimal
{\rm(}and hence not endodualisable{\rm)}.
\end{prop}

Finally we need to consider algebras which have $\bold 2^\ell$ as a
factor, where $\ell \geq 2$ (case~(D)).

\begin{prop} {\rm(}Non-endodualisable finite double Stone
algebras, Case D{\rm)} Let $\bold L = \bold J \times \bold 2^\ell$,
where $\bold J \in \ISP(\bold 4) \setminus \ISP(\bold 3)$ is a
finite double Stone algebra with a non-trivial Boolean core and
$\ell \geq 2$.  Then $\bold L$  is not $3$-endoprimal {\rm(}and so
not endodualisable{\rm)}.
\end{prop}

   We identified firstly various  endodualisable
finite double Stone algebras and then we showed considering in turn
four cases (A)--(D) that there are no other endodualisable finite
double Stone algebras. Here we bring our results together.

\begin{thm} {\rm(}Endodualisability for finite double Stone
algebras, Summary{\rm)} Assume that
 $\bold L = (L;\lor,\land,^{\star} ,^+ ,0,1)$
is a finite proper double Stone algebra with a non-empty core
$K(\bold L)=[a,b]$  {\rm(}$a<b${\rm)}. Then the following are
equivalent:
\begin{itemize}
\item[(1)] $\bold L$ is endodualisable;
\item[(2)] $\bold L$ is endoprimal;
\item[(3)] $\bold L$ is $3$-endoprimal;
\item[(4)] $\bold 5$ is a retract of $\bold L$;
\item[(5)] the core $K(\bold L)$ is a non-Boolean lattice.
\end{itemize}
\end{thm}

For proper double Stone algebras with empty core we have the
following theorem.

\begin{thm} Let $\bold L =  (L;\lor,\land,^{\star} ,^+
,0,1)$ be a finite proper double Stone algebra with empty core. Then
the following are equivalent:
\begin{itemize}
\item[(1)] $\bold L$ is endodualisable;
\item[(2)] $\bold L$ is endoprimal;
\item[(3)] $\bold L$ is $3$-endoprimal;
\item[(4)] $\bold 5\times \bold 2^2$
is a retract of $\bold L$.
\end{itemize}
\end{thm}

For algebras in $\ISP(\bold 3)$ we have, likewise,  the following
result.

\begin{thm}  Let $\bold L$ belong to the variety $\CR
= \ISP(\bold 3)$ of regular double Stone algebras and assume that
$\bold L$ is not Boolean. Then the following are equivalent:
\begin{itemize}
\item[(1)] $\bold L$ is endodualisable;
\item[(2)] $\bold L$ is endoprimal;
\item[(3)] $\bold L$ is $1$-endoprimal;
\item[(4)] $\bold 3 \times \bold 2^2$ is a retract of $\bold L$.
\end{itemize}
\end{thm}

We  record explicitly the following theorem, which  is a corollary
of our  preceding results.

\begin{cor} A finite double Stone algebra is endoprimal
if and only if it is endodualisable.
\end{cor}

\section{Full versus Strong Problem in the theory of natural dualities}

 Every quasi-variety of the form $\CA
= \ISP(M)$, where $\M$ is a finite lattice-based algebra, has a
natural duality. In the case that $M$ is distributive-lattice based,
it is possible to use the \emph{restricted Priestley duality} and
the natural duality for $\CA$ simultaneously. In tandem, these
dualities can provide an extremely powerful tool for the study
of~$\CA$: see Clark and Davey~\cite[Chapter 7]{CD98}. As well as
being a natural area of application of natural duality theory,
distributive-lattice-based algebras in general, and distributive
lattices in particular, have provided deep insights into the general
theory. Important examples have been Heyting algebras, particularly
the finite Heyting chains, and Kleene algebras; but here we firstly
concentrate on the three-element bounded distributive lattice
\[
\bold 3 = (\{0,d,1\};\lor,\land,0,1),
\]
which was seminal in developments that led to the solution of the
\emph{Full versus Strong Problem}, one of the most tantalizing
problems in the theory of natural dualities.

\subsection{The seminal example of the three-element chain}

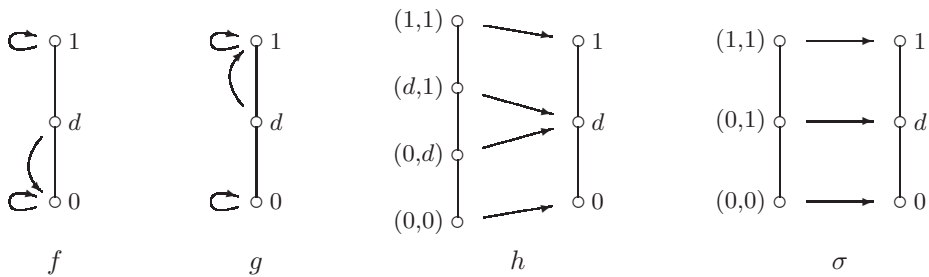
\begin{figure}[ht]
\begin{center}
\begin{picture}(315,100)
\put(0,25){\begin{picture}(0,0)
  \put(0,0){\begin{picture}(0,0)
    \put(0,-25){\makebox(0,0)[b]{\smash{$f$}}}
    \put(0,0){\begin{picture}(0,0)
      \put(0,0){\circle{4}}
      \put(5,0){\makebox(0,0)[l]{\small$0$}}
      \qbezier(0,2)(0,15)(0,28)
      \qbezier(-7,-2)(-17,-5)(-17,0)\qbezier(-7,2)(-17,5)(-17,0)\put(-7,2){\vector(4,-1){0}}
      \end{picture}}
    \put(0,30){\begin{picture}(0,0)
      \put(0,0){\circle{4}}
      \put(5,0){\makebox(0,0)[l]{\small$d$}}
      \qbezier(0,2)(0,15)(0,28)
      \qbezier(-4.881,-5.694)(-15,-17.5)(-4.881,-26.306)\put(-4.881,-26.306){\vector(1,-1){0}}
      \end{picture}}
    \put(0,60){\begin{picture}(0,0)
      \put(0,0){\circle{4}}
      \put(5,0){\makebox(0,0)[l]{\small$1$}}
      \qbezier(-7,-2)(-17,-5)(-17,0)\qbezier(-7,2)(-17,5)(-17,0)\put(-7,2){\vector(4,-1){0}}
      \end{picture}}
    \end{picture}}
  \end{picture}}
\put(75,25){\begin{picture}(0,0)
  \put(0,0){\begin{picture}(0,0)
    \put(0,-25){\makebox(0,0)[b]{\smash{$g$}}}
    \put(0,0){\begin{picture}(0,0)
      \put(0,0){\circle{4}}
      \put(5,0){\makebox(0,0)[l]{\small$0$}}
      \qbezier(0,2)(0,15)(0,28)
      \qbezier(-7,-2)(-17,-5)(-17,0)\qbezier(-7,2)(-17,5)(-17,0)\put(-7,2){\vector(4,-1){0}}
      \end{picture}}
    \put(0,30){\begin{picture}(0,0)
      \put(0,0){\circle{4}}
      \put(5,0){\makebox(0,0)[l]{\small$d$}}
      \qbezier(0,2)(0,15)(0,28)
      \qbezier(-4.881,5.694)(-15,17.5)(-4.881,26.306)\put(-4.881,26.306){\vector(1,1){0}}
      \end{picture}}
    \put(0,60){\begin{picture}(0,0)
      \put(0,0){\circle{4}}
      \put(5,0){\makebox(0,0)[l]{\small$1$}}
      \qbezier(-7,-2)(-17,-5)(-17,0)\qbezier(-7,2)(-17,5)(-17,0)\put(-7,2){\vector(4,-1){0}}
      \end{picture}}
    \end{picture}}
  \end{picture}}
\put(270,25){\begin{picture}(0,0)
  \put(0,0){\begin{picture}(0,0)
    \put(22.5,-25){\makebox(0,0)[b]{\smash{$\sigma$}}}
    \put(0,0){\begin{picture}(0,0)
      \put(0,0){\circle{4}}
      \put(-5,0){\makebox(0,0)[r]{\small$(0{,}0)$}}
      \qbezier(0,2)(0,15)(0,28)
      \end{picture}}
    \put(0,30){\begin{picture}(0,0)
      \put(0,0){\circle{4}}
      \put(-5,0){\makebox(0,0)[r]{\small$(0{,}1)$}}
      \qbezier(0,2)(0,15)(0,28)
      \end{picture}}
    \put(0,60){\begin{picture}(0,0)
      \put(0,0){\circle{4}}
      \put(-5,0){\makebox(0,0)[r]{\small$(1{,}1)$}}
      \end{picture}}
    \end{picture}}
  \put(45,0){\begin{picture}(0,0)
    \put(0,0){\begin{picture}(0,0)
      \put(0,0){\circle{4}}
      \put(5,0){\makebox(0,0)[l]{\small$0$}}
      \qbezier(0,2)(0,15)(0,28)
      \put(-35,0){\vector(1,0){25}}
      \end{picture}}
    \put(0,30){\begin{picture}(0,0)
      \put(0,0){\circle{4}}
      \put(5,0){\makebox(0,0)[l]{\small$d$}}
      \qbezier(0,2)(0,15)(0,28)
      \put(-35,0){\vector(1,0){25}}
      \end{picture}}
    \put(0,60){\begin{picture}(0,0)
      \put(0,0){\circle{4}}
      \put(5,0){\makebox(0,0)[l]{\small$1$}}
      \put(-35,0){\vector(1,0){25}}
      \end{picture}}
    \end{picture}}
  \end{picture}}
\put(150,17.5){\begin{picture}(0,0)
  \put(0,0){\begin{picture}(0,0)
    \put(22.5,-17.5){\makebox(0,0)[b]{\smash{$h$}}}
    \put(0,0){\begin{picture}(0,0)
      \put(0,0){\circle{4}}
      \put(-5,0){\makebox(0,0)[r]{\small$(0{,}0)$}}
      \qbezier(0,2)(0,12.5)(0,23)
      \end{picture}}
    \put(0,25){\begin{picture}(0,0)
      \put(0,0){\circle{4}}
      \put(-5,0){\makebox(0,0)[r]{\small$(0{,}d)$}}
      \qbezier(0,2)(0,12.5)(0,23)
      \end{picture}}
    \put(0,50){\begin{picture}(0,0)
      \put(0,0){\circle{4}}
      \put(-5,0){\makebox(0,0)[r]{\small$(d{,}1)$}}
      \qbezier(0,2)(0,12.5)(0,23)
      \end{picture}}
    \put(0,75){\begin{picture}(0,0)
      \put(0,0){\circle{4}}
      \put(-5,0){\makebox(0,0)[r]{\small$(1{,}1)$}}
      \end{picture}}
    \end{picture}}
  \put(45,7.5){\begin{picture}(0,0)
    \put(0,0){\begin{picture}(0,0)
      \put(0,0){\circle{4}}
      \put(5,0){\makebox(0,0)[l]{\small$0$}}
      \qbezier(0,2)(0,15)(0,28)
      \qbezier(-9.864,-1.644)(-22.5,-3.75)(-35.136,-5.856)
      \put(-9.864,-1.644){\vector(4,1){0}}
      \end{picture}}
    \put(0,30){\begin{picture}(0,0)
      \put(0,0){\circle{4}}
      \put(5,0){\makebox(0,0)[l]{\small$d$}}
      \qbezier(0,2)(0,15)(0,28)
      \qbezier(-9.635,-2.676)(-22.5,-6.25)(-35.365,-9.824)
      \qbezier(-9.635,2.676)(-22.5,6.25)(-35.365,9.824)
      \put(-9.635,-2.676){\vector(4,1){0}}
      \put(-9.635,2.676){\vector(4,-1){0}}
      \end{picture}}
    \put(0,60){\begin{picture}(0,0)
      \put(0,0){\circle{4}}
      \put(5,0){\makebox(0,0)[l]{\small$1$}}
      \qbezier(-9.864,1.644)(-22.5,3.75)(-35.136,5.856)
      \put(-9.864,1.644){\vector(4,-1){0}}
      \end{picture}}
    \end{picture}}
  \end{picture}}
\end{picture}
\end{center}
\caption{The (partial) operations $f$, $g$, $h$ and $\sigma$ on
$3$}\label{fig:fgh}
\end{figure}

For a natural-duality viewpoint, Priestley duality for the class
$\CD$ of bounded distributive lattices is obtained via homsets based
on the two-element chain $\bold 2$ and uses the fact that $\CD =
\ISP(\bold 2)$. By using the fact that $\CD = \ISP(\bold 3)$,
in~\cite{DHP:endo} we introduced the following modified Priestley
duality for $\CD$ as a natural duality based on~$\bold 3$. Let $f,g$
be the non-identity endomorphisms of $3$ (see Figure~\ref{fig:fgh})
and let
\[
\thT = (\{0,d,1\}; f, g, \Tp),
\]
where $\Tp$ is the discrete topology $\Tp$.

Let $\CX = \IScP(\thT)$  be the class of all isomorphic copies of
closed substructures of non-zero powers of $\thT$.

In~\cite{DHP:endo} we showed that such a modified Priestley duality
for $\CD$, in which the order is replaced by endomorphisms, can be
based on any finite non-boolean distributive lattice $\M$. We also
showed that, while the order relation cannot be removed in the
boolean case, it can at least be replaced by any finitary relation
on~$\M$, which itself, like the order on $\bold 2$, forms a
non-boolean lattice.

In~\cite{DH00} we studied the enrichment of $\thT$ given by
\[
\thT_\sigma := (\{0,d,1\}; f, g, \sigma, \Tp),
\]
and in~\cite{DHW05a} we explored deeply the enrichments
$\thT_\sigma$ and
\[
\thT_h := (\{0,d,1\}; f, g, h, \Tp).
\]
(The  binary partial operations $h$ and $\sigma$ are also given in
Figure~\ref{fig:fgh}.) If in the above scheme for the modified
Priestley duality for $\CD$ based on $\bold 3$ the alter ego $\thT$
of $\bold 3$ is replaced with the alter ego $\thT_\sigma$, then not
only the map $e_A: \A \to ED(\A)$ is an isomorphism, for all $\A\in
\CD$, establishing a duality between $\CD = \ISP(\bold 3)$ and
$\CX_\sigma = \IScP(\thT_\sigma)$, but moreover the map
$\varepsilon_X: \X \to DE(\X)$ is an isomorphism, for all $\X\in
\CX_\sigma$, establishing a \emph{full duality} between $\CD$ and
$\CX_\sigma$. In general, such a scheme provides us with a canonical
way of constructing, via hom-functors, a~dual adjunction between a
category of algebras $\CA=\ISP(\M)$, generated by a finite algebra
$\M$, and a category $\CX = \IScP(\MT)$ of structured topological
spaces, generated by the alter ego $\MT$ of the algebra $\M$. (It
should be noted that for some finite algebras $\M$ there is no
choice of alter ego $\MT$ for which the resulting dual adjunction
yields a duality between $\CA$ and $\CX$; for example, the
two-element implication algebra $\I=(\{0,1\};\to)$,
see~\cite[Chapter~10]{CD98}.) If the hom-functors $D, E$ are
restricted to the categories $\CA_{\fin}$ and $\CX_{\fin}$ of finite
members of $\CA$ and $\CX$ only, then the concepts of a
\emph{finite-level} duality, full duality or strong duality are
obtained.

The properties of the modified Priestley dualites for $\CD$ based on
$3$ given by the alter egos $\thT$, $\thT_h$ and $\thT_\sigma$ are
summarized in the following theorem.

\begin{thm}
Let $\thT$, $\thT_h$ and $\thT_\sigma$ be the alter egos of $\bold
3$ defined above.
\begin{itemize}
\item[\rm(i)] $\thT$ yields a duality on $\CD$. {\upshape(Davey, Haviar, Priestley \cite{DHP:endo})}
\item[\rm(ii)]$\thT_h$ yields a full  duality, which is not strong, on the category $\CD_{\fin}$ and yields a duality, which is not full, on the category $\CD$.
{\upshape(Davey, Haviar, Willard \cite{DHW05a})}
\item[\rm(iii)] $\thT_\sigma$ yields a strong duality for $\CD$. {\upshape(Davey, Haviar \cite{DH00})}
\item[\rm(iv)] Every full duality on $\CD$ based on $\bold 3$ is strong. {\upshape(Davey, Haviar, Willard \cite{DHW05a})}
\end{itemize}
\end{thm}

\subsection{Full versus Strong Problem: its local versions and when
full implies strong}

Since the Full versus Strong Problem in its global version had
remained open for the 25 years, we introduced in \cite{DHN} local
versions of this problem that could prove more tractable and
fruitful.

\begin{pr}\label{q:local}
For an arbitrary finite algebra $\M$ in your favourite class $\cat
C$ of algebras, is every full duality based on $\M$ necessarily
strong?
\end{pr}

We also posed the finite-level version of Problem~\ref{q:local}.

\begin{pr}\label{q:local:finite}
For an arbitrary finite algebra $\M$ in your favourite class $\cat
C$ of algebras, is every duality based on $\M$ that is full at the
finite level necessarily strong at the finite level?
\end{pr}

The first solutions to these local versions of the Full versus
Strong Problem were given for full dualities based on the
three-element chain in the variety of bounded distributive lattices
in our paper~\cite{DHW05a} (as shown in the previous subsection).
The answer was shown to be affirmative to Problem~\ref{q:local} and
negative to Problem~\ref{q:local:finite}. In~\cite{DHN} we provided
affirmative answers to Problems~\ref{q:local}
and~\ref{q:local:finite} for full dualities based on an arbitrary
finite algebra in three varieties of algebras: abelian groups,
semilattices (with or without bounds) and relative Stone Heyting
algebras. We also developed some general conditions under which
`full implies strong' that had the potential to add to the list of
solutions. Finally, we answered Problem~\ref{q:local} in the
affirmative for full dualities based on an arbitrary finite lattice
in the variety of bounded distributive lattices.

There is a further, weaker version of Problem~\ref{q:local} that
deserves to be recorded here.

\begin{pr}\label{q:weak}
In your favourite class $\cat C$ of algebras,  is every fully
dualisable finite algebra necessarily strongly dualisable?
\end{pr}

 It should be noted that the
finite-level variant of this question makes no sense since every
finite algebra $\M$ is strongly dualised at the finite level by the
alter ego $\MT= \langle M; H, \T\rangle$, where $H$ consists of all
finitary algebraic partial operation on~$\M$.

We found in~\cite{DHN} several sufficient conditions  for full to
imply strong:

\begin{thm}\label{third:new}
Let\/ $\D$ be a finite algebra, let\/ $\M$ be a finite algebra in
$\CA := \ISP(\D)$ such that\/ $\D$ is a subalgebra of\/~$\M$. Assume
that\/ $\DT = \langle D; G^D,H^D, R^D, \T\rangle$ strongly dualises
$\D$ {\rm[}at the finite level\/{\rm]} and that $D$, each relation
$r\in R^D$, and $\dom(h)$, for all $h\in H^D$,  is an intersection
of equalizers of pairs of algebraic total operations on~$\M$. Then
any alter ego $\MT$ that fully dualises~$\M$ {\rm[}at the finite
level\/{\rm]} strongly dualises~$\M$ {\rm[}at the finite
level\/{\rm]}.
\end{thm}

When $R^D = \varnothing$ there is a particularly satisfying
simplification of this result that involves assumptions on $\D$
only. We say that $\D$ is a \emph{subretract} of $\M$ if $\D$ is a
subalgebra of $\M$ and there is a \emph{retraction} of $\M$ onto
$\D$, that is, a homomorphism $\omega : \M \to \D$ with $\omega
\restriction D = \id_D$.

\begin{thm}\label{third:newer}
Let\/ $\D$ be a finite algebra and let\/ $\CA := \ISP(\D)$. Assume
that\/ $\DT = \langle D; G^D,H^D, \T\rangle$ strongly dualises $\D$
{\rm[}at the finite level\/{\rm]} and that, for all $h\in H^D$, the
set $\dom(h)$ is an intersection of equalizers of pairs of algebraic
total operations on~$\D$. Let $\M$ be a finite algebra in $\CA$ such
that $\D$ is a subretract of\/~$\M$. Then any alter ego $\MT$ that
fully dualises~$\M$ {\rm[}at the finite level\/{\rm]} strongly
dualises~$\M$ {\rm[}at the finite level\/{\rm]}.
\end{thm}

The version of Theorem~\ref{third:newer} that applies when $\DT$ is
a total algebra turned out to be so striking that we stated it as a
separate result:

\begin{thm}\label{second}
Let\/ $\D$ be a finite algebra, let $\CA := \ISP(\D)$ and let\/ $\M$
be a finite algebra in $\CA$ that has\/ $\D$ as a subalgebra. Assume
that\/ $\DT = \langle D; G^D, \T\rangle$ is a total algebra that
strongly dualises $\D$ {\rm[}at the finite level\/{\rm]}. If\/ $\MT$
is an alter ego of\/ $\M$ that fully dualises~$\M$ {\rm[}at the
finite level\/{\rm]}, then $\MT$ strongly dualises~$\M$ {\rm[}at the
finite level\/{\rm]}.
\end{thm}

Also we presented the following special case of
Theorem~\ref{third:new}:

\begin{thm}\label{fourth}
Let\/ $\D$ be a finite algebra. Assume that\/ $\DT = \langle D;
G^D,H^D, R^D, \T\rangle$ strongly dualises $\D$ {\rm[}at the finite
level\/{\rm]} and that each relation $r\in R^D$, and $\dom(h)$, for
all $h\in H^D$, is an intersection of equalizers of pairs of
algebraic total operations on~$\D$. Then any alter ego that fully
dualises~$\D$ {\rm[}at the finite level\/{\rm]}, strongly
dualises~$\D$ {\rm[}at the finite level\/{\rm]}.
\end{thm}

We then applied Theorem~\ref{second} to show that
Questions~\ref{q:local} and~\ref{q:local:finite} have affirmative
answers for arbitrary finite algebras in the varieties of abelian
groups and semilattices.

\smallskip

\noindent {\bf  Abelian groups} Let $\M =\langle M; \cdot, {}^{-1},
1\rangle$ be a finite non-trivial abelian group.  Then there is a
cyclic subgroup $\D$ of $\M$ such that $\D$ is a direct factor of
$\M$ and such that $\D$ and $\M$ generate the same
quasi-variety~$\CA$. Since the total algebra $\DT = \langle D;
\cdot, {}^{-1}, 1, \T\rangle$ yields a strong duality on $\CA$ based
on $\D$ (see~\cite[4.4.2]{CD98}), we may apply Theorem~\ref{second}
to obtain that every alter ego $\MT$ that fully dualises the finite
abelian group $\M$ {\rm[}at the finite level\/{\rm]} also strongly
dualises $\M$ {\rm[}at the finite level\/{\rm]}. Hence the answers
to Questions~\ref{q:local} and~\ref{q:local:finite} in the variety
of abelian groups are always in the affirmative.

\smallskip

\noindent {\bf Semilattices} Let $\D_K=\langle\{0,1\};\lor,K\rangle$
be the two-element semilattice with possible bounds $K\subseteq
\{0,1\}$, let $\cat S_K:=\ISP(\D_K)$ and let $\S$ be a finite
non-trivial semilattice in~$\cat S_K$. We have the following strong
dualities on $\cat S_K:=\ISP(\D_K)$ based on $\D_K$ given by total
algebras.
\begin{itemize}
\item[(i)]
$\DT:=\langle \{0,1\};\lor,0,1,\T\rangle$ yields a strong duality on
$\cat S$ based on the (unbounded) semilattice $\D = \langle
\{0,1\};\lor\rangle$.
\item[(ii)]
$\DT_0=\langle \{0,1\};\lor,0,\T\rangle$ yields a strong duality on
$\cat S_0$ based on the semilattice with zero $\D_0 = \langle
\{0,1\}; \lor, 0\rangle$.
\item[(iii)]
$\DT_1=\langle \{0,1\};\lor,1,\T\rangle$ yields a strong duality on
$\cat S_1$ based on the semilattice with one $\D_1 = \langle
\{0,1\}; \lor, 1\rangle$.
\item[(iv)]
$\DT_{01}=\langle \{0,1\};\lor,\T)$ yields a strong duality on $\cat
S_{01}$ based on the bounded semilattice $\D_{01} = \langle \{0,1\};
\lor,0, 1\rangle$.
\end{itemize}
According to Theorem~\ref{second}, if $\MT$ is an alter ego of $\S$
that fully dualises the finite semilattice $\S$ {\rm[}at the finite
level\/{\rm]}, then $\MT$ also strongly dualises $\M$ {\rm[}at the
finite level\/{\rm]}. So Questions~\ref{q:local}
and~\ref{q:local:finite} have affirmative answers for arbitrary
finite algebras in these varieties of semilattices (with bounds).

\smallskip

\noindent {\bf Bounded distributive lattices}

Let $\CCD$ be the variety of bounded distributive lattices. We
proved in~\cite{DHN} the following theorem, thereby showing that
Question~\ref{q:local} has an affirmative answer for an arbitrary
finite algebra in the variety of bounded distributive lattices.

\begin{thm}\label{DistLat}
Let\/ $\M$ be a finite non-trivial bounded distributive lattice.
If\/ $\MT$ is an alter ego of\/ $\M$ that yields a full duality on
$\CCD$ {\rm(}based on~$\M${\rm)}, then $\MT$ yields a strong duality
on~$\CCD$.
\end{thm}

\subsection{Full versus Strong Problem: related developments and the solution}

The realm of natural dualities that were known to be full but not
strong at the finite level was for some time a very small one,
consisting of a single example. This example, based on the
three-element bounded distributive lattice, was presented in our
paper~\cite{DHW05a}. In our other developments, we extended this
realm to the class of all natural dualities based on an arbitrary
finite non-boolean bounded distributive lattice~\cite{DHNP}.

The results in~\cite{DHW05a} raised new questions and opened up new
research paths within the field of natural dualities. More
precisely, we were led to ask the following questions
(cf.~\cite{DHNP}):

\begin{itemize}
\item[(a)] Could it be that, for a finite algebra that is strongly
dualisable, every full duality on the quasi-variety it generates is
strong?

\item[(b)] What is it about a finite algebra that allows
its full dualities at the finite level to behave so differently from
its full dualities at the infinite level?

\item[(c)] Which finite algebras generate a quasi-variety
for which every duality that is full [at the finite level] is
necessarily strong?

\item[(d)] Which finite algebras have an alter ego that yields a full but
not strong duality at the finite level?
\end{itemize}

As already mentioned, in~\cite{DHN} we proved that, for each finite
abelian group, semilattice and relative-Stone Heyting algebra, every
duality that is full [at the finite level] is strong [at the finite
level], and, for each finite bounded distributive lattice, every
full duality is strong. This provided a partial answer to
Question~(c) and thereby provided examples with which to study
Question~(b). While Question~(a) could be regarded as wild
speculation, it was supported by the limited evidence available to
us. In order to make headway on questions such as these, we felt we
needed a range of examples of finite algebras that possess a full
but not strong duality at the finite level.

In the paper~\cite{DHNP} we addressed Question~(d). More precisely,
we proved the following result:

\begin{thm}\label{DistLat}
Let $\M$ be a finite non-boolean bounded distributive lattice. Then
there is an alter ego $\MT$ of $\M$ such that
\begin{itemize}
\item[(a)] $\MT$ yields a duality  that is not full on the class $\CD$ of all
bounded distributive lattices, yet
\item[(b)]$\MT$ yields a duality that is full but not strong on
the class of finite bounded distributive lattices.
\end{itemize}
\end{thm}

\smallskip

Hence our Problem~\ref{q:local:finite} was shown to have a negative
answer in the variety of bounded distributive lattices by producing
full but not strong dualities at the finite level based on an
arbitrary finite non-boolean lattice.

The authors had hoped to find a conceptual proof of this last
theorem that would indicate possible generalizations beyond
distributive lattices. A natural approach would be to proceed as
follows: let $\M$ be a finite non-boolean bounded distributive
lattice; then $\M$ has the three-element chain $\bold 3$ as a
retract; in~\cite{DHW05a} an alter ego $\twiddle 3$ for $\bold 3$
was given that yields a full but not strong duality at the finite
level; use the retraction from $\M$ onto $\bold 3$ to lift the alter
ego $\twiddle 3$ up to an appropriate alter ego $\MT$ for~$\M$.
 Unfortunately, this turned out to be too simple minded. We
pursued this and many other approaches but to no avail. The
hoped-for conceptual proof eluded us and we were left with the
direct computational proof presented in~\cite{DHNP}. Nevertheless,
our result provided an infinite number of desired examples where
previously there was only one.

Now, at last, we briefly present the much-seeked solution to the
Full versus Strong Problem that was presented by D.\,M.~Clark,
B.\,A.~Davey and R.~Willard~\cite{CDW07}.

Let $\R := (\{0,a,b,1\};t,\lor,\land,0,1)$ be the four-element chain
with $0 < a < b < 1$ enriched with the ternary discriminator
function~$t$. Let $u$ be the partial endomorphism of $\R$ with
domain $\{0, a, 1\}$ given by $u(a)=b$. In~\cite{CDW07} the authors
showed that the algebra $\R$ provides a negative solution to the
Full versus Strong Problem of the theory of natural dualities:

\begin{thm}
The alter ego $\RT_\bot = (\{0,a,b,1\};\graph(u),\Tp)$ yields a full
but not strong duality on $\ISP(\R)$. {\rm (Clark, Davey, Willard
\cite{CDW07})}
\end{thm}

In general, a finite algebra $\M$ admits essentially only one
finite-level strong duality, but can admit many different
finite-level full dualities. The alter egos $\MT$ yielding the
finite-level full dualities for $\ISP_{\fin}(\M)$ form a doubly
algebraic lattice $\CF(M)$ introduced and studied in B.\,A.~Davey,
J.\,G.~Pitkethly and R.~Willard~\cite{DPiW:full}. The following
theorem summarises results in this direction.

\begin{thm}
\
\begin{itemize}
\item[\rm(i)] $|\CF(\M)|=1$ for any finite semilattice, abelian group or relative Stone Heyting algebra $\M$. {\upshape (Davey, Haviar, Niven \cite{DHN})}
\item[\rm(ii)] $\CF(\M)$ is finite for any finite quasi-primal algebra~$\M$; in particular, for the algebra $\R$ defined above, $|\CF(\R)|=17$.  {\upshape (Davey, Pitkethly, Willard \cite{DPiW:full})}
\item[\rm(iii)] The lattice $\CF(\bold 3)$ is non-modular and has size $2^{\aleph_0}$. {\upshape (Davey, Haviar and Pitkethly \cite{DHPi10})}.
\end{itemize}
\end{thm}

%

\end{document}